\documentclass{amsart}
\usepackage{amsmath,amssymb}

\newcommand{\ob}{\ensuremath{\overline}}

\newcommand{\kx}{\ensuremath{k[x_1,x_2,\ldots ,x_r]}}

\newcommand{\gG}{\ensuremath{\mathcal G}}
\newcommand{\pf}{\noindent{\bf Proof.}\quad}
\newcommand{\codim}{\ensuremath{\mathrm{codim}}}

\newcommand{\ann}{\ensuremath{\mathrm{ann}}}

\newtheorem{theorem}{Theorem}
\newtheorem{prop}[theorem]{Proposition}

\newtheorem{defin}[theorem]{Definition}
\newtheorem{counterex}[theorem]{Counter Example}
\newtheorem{cor}[theorem]{Corollary}

\begin{document}

\title{Borel Fixed Initial Ideals of Prime Ideals in Dimension Two}
\thanks{The research of the author was partially supported by the National 
Security Agency.}
\author{Amelia Taylor}
\address{Department of Mathematics\\
Rutgers University\\
Piscataway, NJ 08854, USA}
\email{\tt ataylor@math.rutgers.edu}

\begin{abstract}
We prove that if the initial ideal of a prime ideal is Borel-fixed 
and the dimension of the quotient ring is less than or equal to two, 
then given any non-minimal associated prime ideal of the initial 
ideal it contains another associated prime ideal of dimension one larger.  
\end{abstract}

\maketitle

Let $R = \kx$ be a polynomial ring over a field.  We will say that 
an ideal $I\subseteq R$ has the \emph{saturated chain property} if 
given any non-minimal associated prime ideal $Q$ of $I$ there exists an 
associated prime ideal $P\subseteq Q$ such that 
$\dim(R/P) = \dim(R/Q)+1$.   Therefore given $Q$ an associated 
prime ideal of $I$, there exists a saturated chain of 
associated prime ideals $P_1\subset P_2\subset\cdots P_n = Q$ such 
that $P_1$ is minimal and 
$\dim(R/P_i) = \dim(R/P_{i+1}) + 1$ for $1\leq i\leq n-1$.  

In 1999 Hosten and Thomas~\cite{HT} proved that the 
initial ideal of a toric ideal has this saturated 
chain property.  This type of connectivity does not exist 
in general for initial ideals of prime 
ideals as the following counter-example, provided by Hosten and 
S. Popescu, illustrates.
\begin{counterex}
{\rm Take the toric ideal $\langle xz-a^2,yz-b^2,tz-c^2\rangle$ and 
substitute $z-t$ for $t$, making the ideal 
$\langle xz-a^2,yz-b^2,z^2-tz-c^2\rangle$ which is still 
a prime ideal.  Using the reverse lexicographic order and the variable order 
$x>y>z>t>a>b>c$ the initial ideal of this prime ideal is 
$J=\langle z^2, yz,xz,zb^2,za^2,ytb^2,xtb^2,xta^2\rangle$.  The primary 
decomposition of $J$ is 
$$J=\langle x,y,z\rangle\cap\langle y,t,z\rangle\cap\langle t,z,a^2\rangle\cap\langle z,a^2,b^2\rangle\cap\langle x,y,z^2,a^2,b^2\rangle.$$  
Hence $Ass(R/J)=\{\langle x,y,z\rangle, \langle y,t,z\rangle,\langle t,a,z\rangle,\langle z,a,b\rangle,\langle x,y,z,a,b\rangle\}$.}\end{counterex}  

Hosten and others have since constructed families of prime ideals such 
that the initial ideal of a prime ideal in the family does not have 
the saturated chain property.  However, it was conjectured that lexicographic 
generic initial ideals of homogeneous prime ideals have the saturated chain 
property.  
We prove that if  $P$ is a homogeneous 
prime ideal of $R=\kx$, $\dim(R/P)=2$ and the initial ideal of $P$ is 
Borel-fixed, then the initial ideal of $P$ has the saturated chain property.  
There are many prime ideals with Borel-fixed 
initial ideal since for any prime ideal $P$ the generic initial ideal 
of $P$ is Borel-fixed.  While we restrict the dimension, we do not make any 
assumptions on the monomial order used, nor do we require that the initial 
ideal be generic, only Borel-fixed which is weaker.


We collect some key definitions and properties and then give the main 
theorem.
A \emph{monomial order} 
$\geq$ on a polynomial ring $R=\kx$ over a field 
$k$ is a total order on the monomials in $R$ such that $m\geq 1$ for 
each monomial $m$ in $R$ and if $m_1,m_2,n$ are monomials in $R$ with 
$m_1\geq m_2$ then $nm_1\geq nm_2$.  A monomial order on a polynomial ring 
in several variables generalizes the notion of degree for a polynomial 
ring in one variable.  The \emph{initial term} of an element $f\in R$, 
denoted $in(f)$, is the largest term of $f$ 
with respect to a fixed monomial order.  
Given an ideal $I$ of $R$, the \emph{initial ideal} is defined to be 
$\langle\{in(f):f\in I\}\rangle$, and is denoted $in(I)$.  It should 
be noted that different monomial orders may yield different initial 
ideals, so whenever an initial ideal is referred to, it is assumed a 
monomial order has been fixed.  A \emph{Gr\"obner basis} is a subset 
$\{g_1, \ldots, g_n\}$ of $I$ such that 
$in(I)=\langle in(g_1),\ldots,in(g_n)\rangle$.  

Let $\gG=Gl(r,k)$ be the $r\times r$ invertible matrices over the field 
$k$.  Let $g=(g_{ij})$ be a matrix in $\gG$.  The \emph{Borel subgroup} 
$B=\{(g_{ij})\in\gG\mid g_{ij}=0 \text{ for } j>i\}$ is the subgroup of 
$\gG$ of lower triangular matrices.  We get an automorphism of $R$ from 
$g$ by describing how $g$ acts on the variables.  Define 
$g(x_i)=\sum_{j=1}^r g_{ij}x_j$.  Notice that if $m$ is a monomial 
$x_1^{j_1}\cdots x_r^{j_r}$ then $g(m)=\prod_l(\sum_k g_{lk}x_k)^{j_l}$.  
Define $g(I)=\{g(f)\mid f\in I\}$ to be the action of $g$ on an ideal $I$ 
in $R$.  An ideal $I$ is said to be \emph{Borel-fixed} if $g(I)=I$ for every 
$g\in B$.  Green~\cite{MG} gives a more convenient way, for our purposes, 
of thinking about Borel-fixed ideals through the following definition and 
proposition.

\begin{defin}\label{def:elem}{\rm \cite[Definition 1.24]{MG}  An elementary 
move $e_k$ for $1\leq k\leq n-1$ is defined by $e_k(x^{J})=x^{\hat{J}}$ where 
$\hat{J}=(j_1,\ldots,j_{k-1},j_{k}+1,j_{k+1}-1,j_{k+2},\ldots,j_r)$ and where 
we adopt the convention that $x^{J}=0$ if some $j_m<0$.   }
\end{defin}  

\begin{prop}\label{prop:borel}
{\rm \cite[Proposition 1.25]{MG}} Let $I$ be a monomial ideal.  
The following are equivalent.  \\
(1) If $x^{J}\in I$, then for every elementary move $e_k(x^J)\in I$;\\
(2) $g(I)=I$ for every $g$ belonging to the Borel subgroup;\\
(3) $in(g(I))=I$ for every $g$ in some open neighborhood of the identity in $B$.
\end{prop}

The main property of Borel-fixed ideals that we want to get from this 
definition is that if $x_jm\in I$, where $m$ is any monomial in $R$, then 
$x_im\in I$ for every $i<j$.  Another key property of Borel-fixed ideals 
involves the structure of the associated primes of $R/I$.  

\begin{cor}\label{cor:assprime}{\rm \cite[Corollary 15.25]{Eis}} If $I$ is 
a Borel-fixed ideal in $R$ and $P$ is an associated prime of $R/I$, then 
$P=(x_1,\ldots,x_j)$ for some $j$.  If $Q=(x_1,\ldots,x_t)$ is a maximal 
associated prime then $x_{t+1},\ldots,x_r$ (in any order) is a maximal 
$(R/I)$-regular sequence in $(x_1,\ldots,x_r)$.    \end{cor}

Let $I$ be the initial ideal of a prime ideal $P\subseteq R=\kx$.  
Assume $<$ is the reverse lexicographic 
order.  A theorem of Bayer and Stillman~\cite[Theorem 2.4]{BS}, gives that 
the image of $x_r$, in $R/I$, is a non-zero divisor in $R/I$ if and only 
if the image of $x_r$, in $R/I$, is a non-zero divisor in $R/P$.  Using 
this and the fact that $P$ is a prime ideal, $x_r$ is a non-zero divisor 
in $R/I$ and hence cannot be in any associated prime of $R/I$.  Hence, if 
$codim (I)=r-2$, that is if $\dim(R/P) = 2$, and $I$ is a Borel-fixed ideal, then 
$Ass(R/I)=\{\langle x_1,\ldots,x_{r-2}\rangle\}$ or
 $Ass(R/I)=\{\langle x_1,\ldots,x_{r-2}\rangle,\langle x_1,\ldots,x_{r-1}\rangle\}$.  
Since the initial ideal of a prime ideal is equidimensional~\cite{KS}, 
a similar assessment of the possible choices for the set of associated prime ideals 
of such an initial ideal establishes that for any monomial order 
$\dim(R/P) = 2$ is the first interesting dimension for this question.



\begin{theorem}
Let $R=\kx$ and $P$ be a homogeneous prime ideal in $R$ such that 
$\dim(R/P)=2$.  Fix a monomial order.  Assume the initial ideal $in(P)$ 
is a Borel-fixed ideal.  Then $in(P)$ has the saturated chain property.
\end{theorem}

\pf  The ring $R$ is regular, so $\dim(R/P)=2$ implies $\codim(P)=r-2$.  Let 
$I=in(P)$, then the codimension of $I$ is also $r-2$.  The prime ideal 
$(x_1,\ldots,x_{r-2})$ is minimal over $I$, since $I$ is 
equidimensional~\cite{KS}, Borel-fixed and codimension $r-2$.  
Therefore $(x_1,\ldots,x_{r-2})$ is in $Ass(R/I)$.  We assume  
$(x_1,\ldots,x_r)\in Ass(R/I)$, as otherwise, by Corollary~\ref{cor:assprime}, 
there is nothing to prove.  We prove that $(x_1,\ldots,x_{r-1})\in Ass(R/I)$.   

Since $(x_1,\ldots,x_r)\in Ass(R/I)$ there exists a monomial 
$z\in R\setminus I$ such that $\ann_{(R/I)}(\ob{z})=(\ob{x_1},\ldots,\ob{x_r})$,
 where $\ob{x}$ denotes the image of $x$ in $R/I$.  All such $z$ are 
in the socle of $R/I$ which is a finite dimensional vector space and hence 
there is a monomial of maximal total degree such that its annihilator, in 
$R/I$, is exactly the maximal ideal $(\ob{x_1},\ldots,\ob{x_r})$.  Choose 
$z$ to be a monomial of maximal total degree in the socle of $R/I$.  Since 
$\ann_{R/I}(\ob{z})=(\ob{x_1},\ldots,\ob{x_r})$ we know 
$z(x_1,\ldots,x_r)\subseteq I$.  Hence there exist $f_1,\ldots,f_r\in P$ 
homogeneous, such that $f_i=x_iz+h_i$ and $x_iz>in(h_i)$ for $1\leq i\leq r$.  
We may assume $in(h_i)\notin I$, since if $in(h_i)\in I$, then there exists 
a homogeneous polynomial $G$ in $P$ such that $G=in(h_i)+g$ and $in(G)=in(h_i)$.  
The difference $f_i-G$ is in $P$ and $in(f_i-G)=x_iz$.  Repeat this process until 
a polynomial is obtained where all terms smaller than $x_iz$ are not in $I$ 
replace $f_i$ with this new 
polynomial.  During this process the leading term  of the polynomial 
did not change.  Also, for 
$i=r-1,r$, $f_i$ must have a second term, since if there is not a second term 
then $x_iz\in P$ implies either $x_i\in P$ or $z\in P$.  If $z\in P$ then 
$z\in I$, a contradiction.  Thus $x_i\in P$, but $i=r-1$ or 
$i=r$ so $(x_1,\ldots,x_{r-1})\subseteq I$ contradicts $\codim(I)=r-2$.  
Hence $f_i$ must have a second term for $i=r,r-1$.

Form the S-polynomial for $f_{r-1}$ and $f_r$, that is 
\begin{align}
x_rf_{r-1}-x_{r-1}f_r=& x_rh_{r-1}-x_{r-1}h_r\in P. \label{s:34}
\end{align}
Since the S-polynomial is in $P$ its initial term is in $I$.  

Suppose all of the monomials that form $x_rh_{r-1}$ and those that 
form $x_{r-1}h_r$ are the same.  Then $x_r$ divides $h_r$ and $x_{r-1}$ 
divides $h_{r-1}$.  Therefore $x_{r-1}(z+h_{r-1}/x_{r-1})\in P$ and 
hence $x_{r-1}\in P$ or $z+h_{r-1}/x_{r-1}\in P$.  The first case implies 
$(x_1,\ldots,x_{r-1})\in I$ and the second case implies $z\in I$ which 
are both contradictions.  Write $h_{r-1}=r_1m_1+\ldots+r_sm_s$ where 
$r_1,\ldots r_s\in k$ and $m_1,\ldots,m_s$ are monomials such that 
$m_1>\ldots>m_s$ and similarly write $h_r=s_1n_1+\ldots+s_tn_t$, 
$s_1,\ldots, s_t\in k$, $n_1>\ldots > n_t$ monomials.  Choose 
$i$ smallest such that either $x_rm_i=x_{r-1}n_i$ and $r_i\neq s_i$ or 
$x_rm_i\neq x_{r-1}n_i$.  We consider two cases. \\  
				
(1):  $x_rm_i=x_{r-1}n_i$ and $r_i\neq s_i$ or $x_rm_i>x_{r-1}n_i$ so 
$x_rm_i$ is the leading monomial of the S-polynomial $S(f_{r-1},f_r)$.\\

(2):  $x_{r-1}n_i >x_rm_i$ so $x_{r-1}n_i$ 
is the leading monomial of $S(f_{r-1},f_r)$.\\

In case (1) $x_rm_i$ is in $I$.  This implies $x_jm_i\in I$ for 
$1\leq j\leq r$ since $I$ is Borel-fixed (Proposition~\ref{prop:borel}).  
Hence $m_i(x_1,\ldots,x_r)\subseteq I$ which implies 
$(\ob{x_1},\ldots,\ob{x_r})=\ann_{R/I}(\ob{m_i})$ and hence $\ob{m_i}$ is 
in the socle of $R/I$.  The total degree of $x_{r-1}z$ is the same as that 
for $m_i$ since $f_{r-1}$ is homogeneous, so the total degree of 
$m_i$ is strictly larger than the total 
degree of $z$.  This contradicts our choice of $z$ in the socle.

In case (2) $x_{r-1}n_i$ is in $I$.  Again, $x_in_i\in I$ for $1\leq i\leq r-1$ 
and hence $(\ob{x_1},\ldots,\ob{x_{r-1}})\subseteq \ann_{R/I}(\ob{n_i})$.  Since $n_i\notin I$, if we have equality then 
$(x_1,\ldots,x_{r-1})\in Ass(R/I)$.  If the containment 
$(\ob{x_1},\ldots,\ob{x_{r-1}})\subset \ann_{R/I}(\ob{n_i})(\ob{x_1},\ldots,\ob{x_{r-1}})$ is strict then $\ob{x_r}^t$ is in $\ann_{R/I}(\ob{n_i})$ for some $t\geq 1$.  
Choose $t$ least such that $\ob{x_r}^t\in \ann_{R/I}(\ob{n_i})$.  Then 
$x_r^tn_i\in I$ and hence $x_jx_r^{t-1}n_i\in(I)$ for $1\leq j\leq r$.  
This implies 
\begin{equation}\label{eq:ann_contra}
(\ob{x_1},\ldots,\ob{x_r})=\ann_{R/I}(\ob{x_r^{t-1}n_i}),
\end{equation} 
 and since $x_r^{t-1}n_i\notin I$ by the minimality of $t$, 
$x_r^{t-1}n_i$ is in the socle of $R/I$.  The total 
degree of $x_r^{t-1}n_i$ is greater than or equal to that of $n_i$ 
and hence strictly greater than the total degree of $z$.  Therefore, 
Equation~\ref{eq:ann_contra} gives a contradiction to our choice of $z$ and 
$\ann_{R/I}(\ob{n_i})=(\ob{x_1},\ldots,\ob{x_{r-1}})$ as desired. 
\hfill $\Box$


\end{document}